\documentclass{amsart}
\usepackage{amsthm,amsmath,amssymb,epic,eepic}
{\newtheorem{theorem}{Theorem}[section]                     
        
   \newtheorem{lemma}[theorem]{Lemma}

}
{\theoremstyle{definition}
   \newtheorem{definition}[theorem]{Definition}
   \newtheorem{remark}[theorem]{Remark}
}

\newcommand{\QQ}{{\mathbb{Q}}}

\newcommand{\PP}{{\mathbb{P}}}
\newcommand{\ZZ}{{\mathbb{Z}}}

\newcommand{\cC}{{\mathcal C}}

\newcommand{\cH}{{\mathcal H}}
\newcommand{\cI}{{\mathcal I}}
\newcommand{\cL}{{\mathcal L}}

\newcommand{\cO}{{\mathcal O}}

\newcommand{\cS}{{\mathcal S}}
\newcommand{\cU}{{\mathcal U}}
\newcommand{\cX}{{\mathcal X}}

\newcommand{\Spec}{\operatorname{Spec}}
\newcommand{\Isom}{\operatorname{Isom}}

\newcommand{\Hom}{{\operatorname{Hom}}}

\newcommand{\Aut}{{\operatorname{Aut}}}

\newcommand{\dar}{\downarrow}

\newcommand{\PGL}{{\operatorname{\mathbf{PGL}}}}
\newcommand{\ocM}{\overline{{\mathcal M}}}
\newcommand{\oM}{\overline{{\mathbf M}}}
\newcommand{\ocH}{\overline{{\mathcal H}}}
\newcommand{\curveuparrow}{%
\begin{picture}(4,10)%
\put(0,10){\vector(-1,2){0}}%
\qbezier(0,0)(4,5)(0,10)
\end{picture}}
\begin{document}
\title{Stable maps and Hurwitz schemes in mixed characteristic}
\author[D. Abramovich]{Dan Abramovich}
\thanks{Partially supported by NSF grant DMS-9700520 and by an Alfred
P. Sloan research fellowship}  
\address{Department of Mathematics\\ Boston University\\ 111 Cummington\\
Boston, MA 02215\\ USA} 
\email{abrmovic@math.bu.edu}
\author[F. Oort]{Frans Oort}
\address{Mathematisch Instituut \\
University of Utrecht 	\\
Budapestlaan 6 		\\
NL - 3508 TA UTRECHT 	\\
The Netherlands 	}
\email{oort@math.uu.nl}
\date{\today}

\maketitle

\section{The problem} The Hurwitz scheme was originally conceived as a
parameter space for simply branched covers of the projective line. A variant of
this is a parameter spaec for  simply branched covers of the projective line,
up to automorphisms of $\PP^1$ - the so called unparametrized Hurwitz
scheme. Other variants involve fixing branching types which are not necessarily
simple (see \cite{ddh}). 
A rigorous algebraic definition of the Hurwitz scheme in characteristic 0 was
given by Fulton in \cite{fulton} - the unparametrized version is its quotient
by the 
$\PGL(2)$ action on $\PP^1$. A natural compactification of the unparametrized
Hurwitz scheme  in characteristic 0 was given by Harris and Mumford \cite{hm} -
the space of admissible  covers.  A rigorous treatment using logarithmic
structures can be found in \cite{mochizuki}, who puts some order in the zoo of
variants one can think of (parametrized vs. unparametrized, ordered branch
points vs. unordered  branch points, stack vs. coarse moduli scheme etc.). A
new treatment using twisted principal bundles and stable maps into stacks can
be found in \cite{av,acv}. In \cite{p} the space of admissible covers is
identified as a closed subscheme in the space of stable maps into
$\ocM_{0,n+1}$. 

We know of no existing treatment of Hurwitz schemes in positive or mixed
characteristic when the degree exceeds the characteristic of the residue
fields. In this note we follow Pandharipande's idea in  \cite{p} and define the
compactified Hurwitz scheme as a subscheme of  the space of stable maps into
$\ocM_{0,n+1}$. Variants with target curves of higher genus are defined as
well. 

In order to follow this idea we have to treat moduli of stable maps in mixed
characteristic. This is implicit in \cite{bm} but has not been available in the
litterature in this generality. We show (Theorem \ref{Th:stable-maps}) that for
any base scheme 
$S$, and any integers $g, n, r$ and $d$  there exists an Artin
algebraic stack with finite stabilizers, denoted $\ocM_{g,n}(\PP^r_S,d)$, which
is proper over $S$, parametrizing
stable maps of $n$-pointed curves of genus $g$ into $\PP^r$ over $S$. This
stack admits a projective coarse moduli scheme $\oM_{g,n}(\PP^r_S,d)$. One
immediately derives the existence of a similar stack $\ocM_{g,n}(X,d)$
parametrizing stable maps of degree $d$ into a projective $S$-scheme of finite
presentation $X\subset \PP^r_S$. 

Denote by $\cC_{0,n} \to \ocM_{0,n}$ the universal family of stable
$n$-pointed curves of genus $0$.  We propose to define the (unparametrized)
Hurwitz stack (with ordered 
simple branch points) to be the closure in  
$\ocM_{g,n}(\cC_{0,n}/{\ocM_{0,n}},d)$ of the locus of Hurwitz covers over
$\QQ$: one identifies an 
admissible cover $C \to D$ with ordered simple branchings $P_1,\ldots,P_n$
as a stable map from the $n$-pointed curve $(C,Q_1,\ldots Q_n)$ to the stable
$n$-pointed rational curve $(D,P_1,\ldots P_n)$, where $Q_i$ are the
ramification points. Here the pointed curve $(D,P_1,\ldots P_n)$ is
identified uniquly as a fiber of $\cC_{0,n} \to \ocM_{0,n}$.

There are some interesting things which happen in this construction. For the
sake of comparison, recall that in characteristic $0$
\begin{enumerate} 
\item the compactified Hurwitz scheme is {\bf irreducible;}
\item it maps {\bf finitely} to $\ocM_{0,n}$;
\item it parametrizes only {\bf finite morphisms} $C \to D$;
\item all the morphisms are {\bf simply branched away from the nodes;} and
\item the generic point corresponds to a morphism where $C$ is {\bf
irreducible.} 
\end{enumerate}

We study in some detail the reduction in characteristic 2 of the compactified
stack of 
double covers  $C \to D$ of $\PP^1$ branched in 4 points $0,1,\infty,
\lambda$. It has the
following features:
\begin{enumerate} 
\item this stack is {\bf reducible}: it has four irreducible components; 
\item it {\bf does not map finitely} to $\ocM_{0,n}$; specificatlly, one
component is entirely devoted to maps with $j$ invariant $0$ and arbitrary
$\lambda$ invariant, whereas the other components parametrize maps with
arbitrary $j$ invariant and $\lambda$ invariant $0,1$ or $\infty$;
\item  all geometric points correspond to maps which are {\bf not finite;}
\item  all geometric points correspond to maps whose finite part is
{\bf inseparable}; 
\item the generic points correspond to  morphisms where $C_1$ is {\bf
reducible.} 
\end{enumerate}

\section{Stable maps}
\subsection{Definitions.}
Fix a base scheme $S$, and integers $g,n,r$ and $d$.
\begin{definition}
Let $T$ be an $S$ scheme. A {\em stable, $n$-pointed map $(C\to T, s_i, f)$ of
genus $g$ and 
degree $d$ into $\PP^r_S$ over $T$} consists of a diagram
$$
  \begin{array}{ccc} C & \stackrel{f}{\to} & \PP^r \\
	\pi\dar\curveuparrow s_i&     & \\
		    T & & 
  \end{array}
$$
such that
\begin{enumerate} 
\item the morphism $\pi:C \to T$ is a projective flat family of curves;
\item all the geometric fibers of $C \to T$ are reduced with at most nodes as
singularities;
\item the sheaf $\pi_*\omega_{C/T}$ is locally free of rank $g$;
\item the  $n$ morphisms $s_i:T \to C$ are sections of $\pi$ which are disjoint
and land in the smooth locus of $\pi$;
\item the degree of $ f^*\cO(1)$ on geometric fibers of $C \to T$ is $d$; and
\item\label{Cond:Aut} the group scheme $\Aut_{\PP^r_T}(f:C\to \PP^r)$ is finite
over $T$. 
\end{enumerate}
\end{definition}

\begin{remark}
Denote by $S_i$ the image of $s_i$. The {\em stability} condition
\ref{Cond:Aut} on the 
automorphisms 
is equivalent to the  
condition that 
$\omega_{C/T}(\sum S_i)\otimes f^*\cO(3)$ be ample (see \cite{k},
\cite{fp}). This is also equivalent to 
the usual ``three point condition'' for components of $C$ mapping via $f$ to a
point. 
\end{remark}

 Stable maps form category:
\begin{definition}
Given two $S$ schemes $T, T'$ and stable maps $(C\to T, s_i, f)$ and $(C'\to
T', s'_i, f')$, a {\em morphism} $\alpha$ of stable maps is a commutative
diagram 
$$
  \begin{array}{ccc} C & \stackrel{\alpha_C}{\to} & C' \\
		\dar &	&	\dar \\
		T & \stackrel{\alpha_T}{\to} & T'
  \end{array}
$$
inducing an {\em isomorphism} $C \to C'\times_TT'$, compatible with $f$ and
$s_i$, namely: $\alpha_C\circ s_i = s_i' \circ \alpha_T$, and $f = f'\circ
\alpha_C$.  
\end{definition}

Denote the category of stable, $n$-pointed maps of genus $g$ and
degree $d$ into $\PP^r_S$ by $\ocM_{g,n}(\PP^r_S,d)$.

\subsection{Embedded stable maps and smooth parametrization of
$\ocM_{g,n}(\PP^r_S,d)$.} We follow the axioms of \cite{artin} for 
algebraic stacks. One can easily show directly that
$\ocM_{g,n}(\PP^r_S,d)$ is a stack, and that the ${\cI}som$ functor is
representable. However, since we want to show that it is an {\em algebraic}
stack, it 
is convenient to exhibit it as a quotient stack. This is done in a standard
manner using embeddings in projective space. 

In the rest of this section the discrete parameters $g,n,r,d$ of stable maps
are fixed. Choose two positive integers $M,D$. 

\begin{definition} An {\em embedded stable map} $(C\subset \PP^M_T, s_i, f)$
of  
embedding degree $D$ is a subscheme $C \subset \PP^M_T$, flat over $T$ of
degree $D$, and 
a stable map $(C\to T, s_i, f:C \to \PP^r)$. The embedded stable map is {\em
nondegenerate} 
if the embedding of all the geometric fibers of  $C\to T$ are nondegenerate
(nemely not contained in proper linear projective subspaces).
\end{definition}

 Embedded stable maps form a category: a morphism $(C\subset \PP^M_T, s_i, f)
 \to (C'\subset \PP^M_{T'}, s'_i, f')$ is {\em the pullback  morphism}
 $$C = \PP^M_T\mathop{\times}\limits_{\PP^M_{T'}}C' \to T$$ inducing a morphism
of stable maps $(C\to T, s_i, f)
 \to (C'\to T', s'_i, f')$.

\begin{lemma} Given a fixed projective space $\PP^M$ and an integer $D$, the
category $\ocM_{g,n}((\PP^r\times\PP^M)_S,(d,D))_{\subset \PP^M}^{nd}$ of  
non-degenerately embedded stable maps $(C\subset \PP^M_T, s_i, f)$ with
embedding degree $D$  is a stack
representable by a quasi-projective $S$-scheme $H$. There is 
a natural action of  $\PGL(M+1)$ on this scheme, by translation of the
embedding $C \subset \PP^M$. 
\end{lemma}

{\bf Proof.} First consider the Hibert scheme $\cH_{P,M}$ of subschemes of
$\PP^M$ having Hilbert polynomial $P(T)=D\cdot T-g+1$. Let $\cU_{P,M}\to
\cH_{P,M}$ be the universal 
family. There is a natural action of $\PGL(M+1)$ on $\cU_{P,M}\to \cH_{P,M}$.

We add the data of points: 
there is a closed subscheme $H_1 \subset \cH_{P,M} \times (\PP^M)^n$
parametrizing collections $(X,P_1,\ldots P_n)$ where $X \subset \PP^M$ has
Hilbert 
polynomial $P$ and the points $P_i \in \PP^M$ satisfy $P_i\in X$. There is a
natural diagonal action of  $\PGL(M+1)$ on $H_1$.

A small deformation of a nodal curve is nodal; similarly, a small deformation
of a pointed curve with distinct points has distinct points, and a small
deformation of a non-degenerate embedding is non-degenerate. Thus, there is an
open  
quasi-projective subscheme $H_2 \subset H_1$, stable under the  action of
$\PGL(M+1)$, parametrizing embedded pointed curves $(X,P_1,\ldots P_n)$  
such that $X$ is a nodal curve, which spans $\PP^M$, and the points
$P_i$ are 
distinct. Denote the 
pullback of $\cU_{P,M}\to \cH_{P,M}$ by $U_2\to H_2$.

We now add the data of the map $f$:
By \cite{Grothendieck}, there is a quasi projective scheme $\Hom_{H_2}(U_2,
\PP^r)$, of finite presentation, parametrizing 
morphisms $f:C \to  \PP^r$ where $C$ is a fiber of $U_2\to H_2$ and $f^*\cO(1)$
has degree $d$. The action of  $\PGL(M+1)$ lifts by composing with $f$.

Finally, we impose the stability condition: there is an open subscheme
$H\subset H_3$ over which the sheaf 
$\omega_{U/H_3}(\sum S_i)\otimes f^*\cO(3)$ is ample. The scheme $H$
satisfies the requirements of the lemma.
\qed

We will now choose the integers $M$ and $D$:

Fixing an integer $\nu\geq 3$, it was shown in \cite{fp} that the invertible
sheaf $$\cL^\nu
=\left(\omega_{C/T}(\sum S_i)\otimes  
f^*\cO(3)\right)^\nu$$
 is very ample and has no higher cohomology along the
fibers. 
Let $\dim H^0(C_t,\cL^\nu) = M+1 $ and $\deg_{C_t}{\mathcal
L}^\nu = D$. The lemma above provides us with a scheme $H$ parametrizing
non-degenerately embedded stable maps with embeding degree $D$ in $\PP^M$.
By \cite{Mumford}, 
 there is a closed subscheme $V$ where the 
embedding sheaf coincides along the fibers with $\cL^\nu$. There is a
natural $\PGL(M+1)$ action on $V$. 

We claim that the stack quotient $V/\PGL(M+1)$ is equivalent to
$\ocM_{g,n}(\PP^r_S,d)$. Indeed, there is an obvious forgetful functor
$V/\PGL(M+1) \to \ocM_{g,n}(\PP^r_S,d)$, by leaving out the embedding. On the
other hand, given a stable map $(\pi:C\to T, s_i, f)$, the projective frame
bundle $P\to T$ associated with the locally free sheaf $\pi_*\cL^\nu$ is a
principal 
$\PGL(M+1)$-bundle over which $\PP(\pi_*\cL^\nu)$ has a canonical
trivialization, therefore there is a $\PGL(M+1)$-equivariant morphism $P \to
V$,  giving rise to a morphism $T  \to V/\PGL(M+1)$. This induces a functor
$\ocM_{g,n}(\PP^r_S,d) \to V/\PGL(M+1)$.  The  two compositions of these
functors are 
easily seen to be equivalent to the identity.

Thus $\ocM_{g,n}(\PP^r_S,d)$ is an Artin algebraic stack.

\subsection{Properness}
Since $\ocM_{g,n}(\PP^r_S,d)$ is a quotient stack it is already of finite
presentation. To check that $\ocM_{g,n}(\PP^r_S,d) \to S$ is proper, we need to
verify the valuative criteria for separation and properness. 

Let $T$ be an $S$ scheme which is the spectrum of a discrete valuation ring
$R$, 
with special point $s$ and generic point $\eta$. Assume we are given two stable
maps  $(\pi:C\to T, s_i, f)$ and  $(\pi_1:C^1\to T, s^1_i, f^1)$  over the same
scheme 
$T$, 
and an isomorphism along the generic fiber $\alpha_\eta:(C\to T, s_i, f)_\eta
\to (C^1\to T, s^1_i, 
f^1)_\eta$. We wish to extend this isomorphism over $T$. 

The closure $C'$ of the graph
of $\alpha$ is proper and flat over 
$T$. It admits maps $r:C'\to C$ and $r^1:C'\to C^1$.  Considering the morphisms
$f\circ r$ and $f^1 \circ r^1$, we see that their graphs are the closure of the
graph of the same morphism on the generic fiber; therefore $f\circ r = f^1
\circ r^1$. Call this map 
$f'$. Similarly, the sections  
$s_i$ and $s_i^1$ lift to sections $s_i':T\to C'$. We may replace $C'$ by a
desingularization on which the images of the sections $S_i'$ and the central
fiber $C'_s$ together form a divisor of normal crossings. 

We claim  that $C$ is the image
of $C'$ under the relative linar series of $(\pi\circ r)_*\left
( (\omega_{C'/T}(\sum 
S'_i)\otimes  {f'}^*\cO(3))^\nu\right)$. Indeed, $\omega_{C'/T}(\sum
S'_i)\otimes  {f'}^*\cO(3))^\nu = \omega_{C/T}(\sum
S_i)\otimes  {f}^*\cO(3))^\nu(E)$ where $E$ is an effective $r$-exceptional
divisor, and therefore we have 
$$(\pi\circ r)_*\left
( (\omega_{C'/T}(\sum 
S'_i)\otimes  {f'}^*\cO(3))^\nu\right) = (\pi)_*\left
( (\omega_{C/T}(\sum 
S_i)\otimes  {f}^*\cO(3))^\nu\right)$$
Since the same holds for $C^1$, we have an
isomorphism $C \to C^1$. It follows that the $\Isom$ scheme is proper. Thus
$\ocM_{g,n}(\PP^r_S,d) \to S$ is separated. Since
stable maps have finite automorphisms, the $\Isom$ scheme is in fact finite.

Continuing to work with $T$ as above, let $(C_\eta\to \eta, (s_i)_\eta,
f_\eta)$ be a stable map. We wish to extend it over $T$, at least
after a finite base change $T'\to T$.  

There always exists a projective extension $C_0 \to T$ of $C_\eta$, by taking
the closure in some projective embedding of $C_\eta$. Replacing
$C_0$ by the closure of the graph of $f_\eta$ in $C_0\times\PP^r$ we may assume
$f_\eta$ extends to $C_0$. 
%By Abhyankar's theorem (see \cite{Lipman-res}) we
%may assume $C_0$ non-singular. 
The sections  $(s_i)_\eta$
are extended by taking the closure of $(S_i)_\eta$.
%as disjoint sections, by making their union with the central fiber
%into a divisor of normal crossings.   
%Our next goal is to repalce $C_0$ by a family of nodal curves.

\begin{lemma}
There exists a finite extension of discrete valuation rings $R \subset R'$
giving rise to a finite base change $T'\to T$, a diagram
$$
	\begin{array}{ccc}
		C_\eta\mathop{\times}\limits_T T' & \subset & C' \\
		\dar				&	&\dar \\
		\{\eta'\}			& \subset & T'	
	\end{array}
$$
and an extension $f':C'\to \PP^r$ of $f_\eta$ and $s_i:T'\to C'$ of
$(s_i)_\eta$, such that $C'\to T'$ is a family of pointed nodal curves. 
\end{lemma}

%We provide two proofs of this lemma.

{\bf Proof.} First assume $C_\eta$ is smooth. The results of \cite{dj2} say in
particular that there exists a 
scheme $T'$ as in the lemma, and a diagram
$$
	\begin{array}{ccc}
		C' & \to  & C_0 \\
		\dar & &   \dar \\
		T' & \to  & T
	\end{array}
$$
such that $C' \to C_0\times_TT'$ is birational and $C'\to T'$ is a family of
pointed nodal curves. This immediately 
proves the lemma in this case. 

For the general case, first take a base change so that the irreducible
components of  $C_\eta$ are absolutely irreducible, and all the nodes are
split. Apply the argument above to the normalization of each irreducible
component of  $C_\eta$, adding extra sections for the points above the nodes.
Finally, glue the resulting schemes together along the sections as in
\cite{bm}. The glued 
object 
is an Artin algebraic space.

Now $(C'\to T', s_i, f)$ is not stable - we need to contract the so called
``components 
to be contracted'' of \cite{bm}. One can apply a modified verision of Knudsen's
contraction technique 
directly to carry that through. Alternatively, after a suitable base change one
can add enough sections meeting the components whose image in $\PP^r$ is
nontrivial, and then Knudsen's contraction technique for stable pointed curves
applies without modifications.

Thus $\ocM_{g,n}(\PP^r_S,d)\to S$ is proper.

\begin{remark} A proof of De Jong's result which we used above in this case can
be summarized as 
follows: one first picks a projective completion $X \to T$ of $C_\eta$ with a
morphism $X \to \PP^r$. After a 
base change, one adds enough sections so that every component 
of every fiber has at least three sections in the non-singular locus of its
reduction.  After a further base change, one has a stable pointed curve $C'\to
T'$ and  a rational map $C' \to X$. De Jong's ``three point lemma'' (which is
relatively straightforward in the case of one-dimensional base) implies
that this is a morphism.

As noted in \cite{fibered}, De Jong's alteration theorems can be deduced from
the existence of the moduli stack of stable maps. 
\end{remark}

\subsection{Projectivity of the coarse moduli scheme.} By \cite{keel-mori}
there exists a coarse moduli morphism
$\ocM_{g,n}(\PP^r_S,d)\to\oM_{g,n}(\PP^r_S,d)$  such that
$\oM_{g,n}(\PP^r_S,d)\to S$ is a proper Artin algebraic space. We 
endow the category $\ocM_{g,n}(\PP^r_S,d)$ with a canonical semipositive
polarization (\cite{kollar}, 2.3 and 2.4). This implies that
$\oM_{g,n}(\PP^r_S,d)\to S$ is projective by \cite{kollar}, Theorem 2.6. Note
that the assumption on tame automorphisms in \cite{kollar} was used only to
guarantee the existence of the algebraic space $\oM_{g,n}(\PP^r_S,d)$; however,
by \cite{keel-mori} it is enough that the automorphism schemes  are finite.

In order to define such a canonical polarization, one could modify the proof of
\cite{kollar}, Proposition 4.7 to apply for the 
sheaf $\cL^\nu$ we constructed above for
stable maps. Alternatively, we can reduce to the case treated by Koll\'ar as
follows: 

First note that since $\ocM_{g,n}(\PP^r_S,d)$ is proper, the inseparable degree
of any such stable map $(C \to  T, s_i, f)$ is bounded from above by some
integer $B>1$.  

 We pick the canonical polarization $\cL^\nu$  as above with $\nu= 
B$. We claim that this is a semipositive polarization.  Let $(C \to 
T, s_i, f)$ be a stable map over a projective curve $T$.  Choose a general
hypersurface  $H$ in $\PP^r$ of degree $3B$. In particular we may assume that
the  pullback $f^*H$ is finite over 
$T$ of inseparable degree $<B$. After a 
base change we may assume $f^*H = \sum b_i \Sigma_i$ where $\Sigma_i$ are
sections of $C \to T$. By assumption we have $b_i <B$. Therefore
$(\omega_{C/T}(\sum \cdot S_i)\otimes f^*\cO(3))^\nu =  
  \omega_{C/T}^\nu(\sum \nu \cdot S_i + \sum b_i  \Sigma_i)$ satisfies the
assumptions of   \cite{kollar}, Proposition 4.7, and therefore the polarization
is semipositive. 

We have thus proven:
\begin{theorem}\label{Th:stable-maps}
The category $\ocM_{g,n}(\PP^r_S,d)\to S$ is aproper  Artin
algebraic stack with finite stabilizers, admitting a projective coarse moduli
scheme $\oM_{g,n}(\PP^r_S,d)\to S$. 
\end{theorem} 

\subsection{Stable maps into a projective scheme.} Let $X \subset \PP^r_S$ be a
projective scheme of finite presentation over $S$. Given a stable map $(C\to T,
s_i, f:C \to  \PP^r)$, it follows from \cite{Grothendieck} that  there is a
closed subscheme $T'\subset T$ of finite presentation where the maps land
inside $X$. It follows that there exists a proper stack of finite presentation
$\ocM_{g,n}(X/S,d)\to S$ admitting a projective coarse moduli scheme. As in
\cite{fp}, given an element  $\beta\in A_1(X)$, there is an open and closed
substack $\ocM_{g,n}(X/S,\beta)\subset \ocM_{g,n}(X/S,d)$ where the maps have
image class $\beta$.

\subsection{Canonical maps.} Assume that $2g-2+m>0$ and $n\geq m$. As noted in
\cite{fp} and 
\cite{bm} 
there is a morphism $\ocM_{g,n}(\PP^r_S,d)\to \ocM_{g,m}$ by ``forgetting the
map and the last markings, and contracting the components to be
contracted''. More generally, given an $S$-morphism  $X
\to Y$ of projective schemes over
$S$, there is a morphism  $\ocM_{g,n}(X/S,\beta) \to
\ocM_{g,n}(Y/S,f_*\beta)$. The condition  $2g-2+m>0$ can be removed as long as
$f_*\beta\neq 0$.

\subsection{Stable maps over a base stack.} Let $S_1\to S_2$ be an \'etale
morphism. 
 It is easy to see that  $\ocM_{g,n}(\PP^r_{S_1},d) =
 \ocM_{g,n}(\PP^r_{S_2},d)\times_SS_1$.  Thus
relative stable maps are well-behaved in the \'etale site. We now apply this to
Deligne-Mumford stacks.

Let $\cS$  be a Deligne-Mumford 
algebraic 
stack. Fix an
\'etale parametrization $S \to \cS$, and let $R= S \times_\cS S$. The remarks
above show that the representable morphism $\ocM_{g,n}(\PP^r_{R},d) \to 
\left(\ocM_{g,n}(\PP^r_{S},d)\right)^2$ is an \'etale groupoid in Artin
stacks. Denote its quotient by  $\ocM_{g,n}(\PP^r_{\cS},d)\to \cS$. One checks
easily 
that this is a proper Artin algebraic stack admitting a relatively coarse
projective moduli stack $\oM_{g,n}(\PP^r_{\cS},d)\to \cS$. More concretely, it
is the quotient of the following groupoid using non-degenerately embedded
stable maps: 
$$\ocM_{g,n}((\PP^r\times \PP^M)_{R},(d,D))_{\subset \PP^M}^{nd}\times\PGL(M+1)
\to  \left(\ocM_{g,n}((\PP^r\times \PP^M)_{S},(d,D))_{\subset
\PP^M}^{nd}\right)^2.$$ 
It is equivalent to
the category of 
stable maps $(C \to T, s_i, f:C \to \PP^r)$ where $T$ runs over schemes over
$\cS$. Note that, since we work with stable maps  {\em relative to a projective
morphism,} there is no need to invoke the setup of \cite{av}.

Asimilar construction works for $\cX \subset \PP^r_{\cS}$, a projective
substack of finite presentation.

\section{Complete Hurwitz stacks} Fix integers  $d\geq 1$ and $n\geq
3$. Let 
$T$ be a 
scheme of pure characteristic 
0. Consider a diagram 
$$\begin{array}{ccccc} C & & \longrightarrow & & D \\
 & \searrow & & \swarrow& \\
 & & T & & \end{array} $$
where $C \to  T$ is a smooth projective connected curve 
admitting $n$ disjoint sections $s_i:T \to C$; where $D \to  T$ is a smooth
projective connected curve of genus $h$, admitting $n$ 
disjoint sections $\sigma_i:T \to D$; and $f:C \to C$ is a Hurwitz cover
simply ramified precisely along $s_i$ with corresponding branch points
$\sigma_i$. The genus $g$ of $C$ is given by the Hurwitz formula $2g-2 =
n+(2h-2)d$. It is shown in \cite{mochizuki} that such diagrams, with morphisms
given 
by pullback diagrams, form a Deligne-Mumford stack $\cH_{d,h,n}$. 

 As noted in \cite{p}, the pointed curve $(D, \sigma_i)$ is a stable pointed
curve of genus 
$h$, and therefore it is the pullback of $\cC_{h,n}$ along a unique morphism
$T \to \ocM_{h,n}$. Composing with the pullback map, we obtain a well defined
morphism $f: C \to \cC_{h,n}$. We thus have a stable $n$-pointed map $(C \to T,
s_i,f: C \to \cC_{h,n})$, of genus $g$, and image class $\beta=d\cdot F$ ($F$
being 
the class of a fiber) over $S=\ocM_{h,n}$. We therefore obtained a functor  
 $\cH_{d,h,n} \to \ocM_{g,n}(\cC_{h,n}/\ocM_{h,n},\beta)$, which is clearly an
embedding of stacks. Denote by  $\ocH_{d,h,n}$ the closure. In case $h=0$ we
write $\ocH_{d,n}$ for $\ocH_{d,0,n}$.

\begin{definition}
The stack  $\ocH_{d,h,n}$ is called {\em the complete Hurwitz stack}.
\end{definition}

There are morphisms $\ocH_{d,h,n} \to \ocM_{h,n}$  (the base curve morphism)
and  
$\ocH_{d,h,n} \to \ocM_{g,n}$ (forgetting the map and contracting the
extraneous components). 

\section{Double covers of $\PP^1$ branched over $4$ points.}
Consider $\ocH_{2,4}$ - the complete Hurwitz stack of double covers of $\PP^1$
branched over $4$ points. Over $\Spec\ZZ[1/2]$ it is a $\ZZ/2$-gerbe over
$\ocM_{0,4}=\PP^1$, equivalent to the stack of elliptic curves with level 2 
structure. Here we are interested in its stucture in characteristic 2. 

We will make use of the morphisms $\lambda:\ocH_{2,4}\to \ocM_{0,4}$ and
$j:\ocH_{2,4}\to\ocM_{1,1}$. In addition, if $(C \to T, s_i, f: C \to D)$ is a
Complete Hurwitz  map   whose generic fiber is in  $\cH_{d,n}$, there is a
natural involution $\sigma$ acting 
 on the generic fiber, which extends to $C$, since $C$ is normal. Consider
the associated stable 1-pointed curve $\overline{C} \to T$ of genus 1 obtained
by 
contracting the extraneous components. The induced involution on the
smooth locus of $\overline{C} \to T$ coincides with
the involution $(-1)$ of the group law, therefore its fixed point locus is the
2-torsion group-scheme.

We also have a natural action of the symemtric group of 4 letters on
$\ocH_{2,4}$, permuting the markings.

Finally, it is easy to check that for any object $(C \to T, s_i,
f: C \to D)$ of the complete Hurwitz stack, the source curve $C$ is a stable
4-pointed curve (and in particular $\ocH_{2,4}$ is a Deligne-Mumford
stack). This follows  
since for any component $C_0$ of $C$ mapping finitely to a component $D_0$ of
$D$, and any marked 
point $P$ of $D_0$, there is a unique point $Q$ of $C_0$ which is either marked
or possibly a node.

Let $(C_s , s_i, f: C_s \to D_s)$ be a geometric point of the Complete Hurwitz 
stack in characteristic 2. We can choose a discrete valuation ring $R$  
 with spectrum $T$, with generic point of characteristic 0 and residue
characteristic $2$, and a stable map $(C\to T , s_i, f: C \to D)$ specializing
to the given one. Denote by $j_s$  and $\lambda_s$ the $j$ invariant and the
$\lambda$ invariant of the special fiber. 

{\sc Case 1:}  $\lambda_s \neq 0, 1, \infty$. 

The curve $C_\eta$
admits a Legendre equation $y^2 = x(x-1)(x-\lambda)$.  This equation has a
unique singular point in characteristic 2 corresponding to $x^2 = \lambda$.  
It is easy to see that this is the unique point invariant under the symmetries
of $(0,1,\infty,\lambda)$: the permutation $0\leftrightarrow
\infty, 1\leftrightarrow \lambda$ is given by the transformation $x \mapsto
\lambda/x$, whose unique fixed 
point is $\sqrt{\lambda}$. The permutation $0\leftrightarrow 1
\infty \leftrightarrow \lambda$ is given by the transformation $ x \mapsto
(x-1)/(\lambda^{-1}x-1)$ with the same fixed point!
 In any case, in the stable map $f_s$ the singular
point  $x^2 = \lambda$ must be blown up (it is not a node!) and for stability
reasons  a
component of genus 1 must be attached.  The reduction of the 2-torsion locus on
$\overline{C}$ consists precisely of  the attaching point, therefore the
elliptic curve is supersingular, namely 
$j_s=0$.

The resulting picture
is as follows: $C_s$ has exactly $2$ components: $C_0$, which has genus 0, and
$C_1$, which has genus 1 and $j$ invariant 0. The
morphism $f_s$ maps $C_1$ to $x = \sqrt{\lambda}$. The component $C_0$
is the normalization of $y^2 = x(x-1)(x-\lambda)$, mapping 2-to-1 to $D_s$,
purely inseparably. It is attached over 
$x = \sqrt{\lambda}, y=\lambda+x$ with the elliptic curve $C_1$.
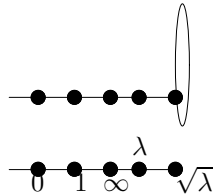
\begin{figure}[h]
$$
\setlength{\unitlength}{0.0005in}
\begingroup\makeatletter\ifx\SetFigFont\undefined
% extract first six characters in \fmtname
\def\x#1#2#3#4#5#6#7\relax{\def\x{#1#2#3#4#5#6}}%
\expandafter\x\fmtname xxxxxx\relax \def\y{splain}%
\ifx\x\y   % LaTeX or SliTeX?
\gdef\SetFigFont#1#2#3{%
  \ifnum #1<17\tiny\else \ifnum #1<20\small\else
  \ifnum #1<24\normalsize\else \ifnum #1<29\large\else
  \ifnum #1<34\Large\else \ifnum #1<41\LARGE\else
     \huge\fi\fi\fi\fi\fi\fi
  \csname #3\endcsname}%
\else
\gdef\SetFigFont#1#2#3{\begingroup
  \count@#1\relax \ifnum 25<\count@\count@25\fi
  \def\x{\endgroup\@setsize\SetFigFont{#2pt}}%
  \expandafter\x
    \csname \romannumeral\the\count@ pt\expandafter\endcsname
    \csname @\romannumeral\the\count@ pt\endcsname
  \csname #3\endcsname}%
\fi
\fi\endgroup
{\renewcommand{\dashlinestretch}{30}
\begin{picture}(2998,2012)(0,-10)
\put(1812,1352){\ellipse{150}{1274}}
\put(1737,1014){\blacken\ellipse{150}{150}}
\put(1737,1014){\ellipse{150}{150}}
\put(312,1014){\blacken\ellipse{150}{150}}
\put(312,1014){\ellipse{150}{150}}
\put(687,1014){\blacken\ellipse{150}{150}}
\put(687,1014){\ellipse{150}{150}}
\put(1062,1014){\blacken\ellipse{150}{150}}
\put(1062,1014){\ellipse{150}{150}}
\put(1362,1014){\blacken\ellipse{150}{150}}
\put(1362,1014){\ellipse{150}{150}}
\put(312,264){\blacken\ellipse{150}{150}}
\put(312,264){\ellipse{150}{150}}
\put(687,264){\blacken\ellipse{150}{150}}
\put(687,264){\ellipse{150}{150}}
\put(1062,264){\blacken\ellipse{150}{150}}
\put(1062,264){\ellipse{150}{150}}
\put(1362,264){\blacken\ellipse{150}{150}}
\put(1362,264){\ellipse{150}{150}}
\put(1737,264){\blacken\ellipse{150}{150}}
\put(1737,264){\ellipse{150}{150}}
\path(12,1014)(1812,1014)
\path(12,264)(1812,264)
\put(237,39){\makebox(0,0)[lb]{\smash{{{\SetFigFont{10}{14.4}{rm}$0$}}}}}
\put(687,39){\makebox(0,0)[lb]{\smash{{{\SetFigFont{10}{14.4}{rm}$1$}}}}}
\put(987,39){\makebox(0,0)[lb]{\smash{{{\SetFigFont{10}{14.4}{rm}$\infty$}}}}}
\put(1287,414){\makebox(0,0)[lb]{\smash{{{\SetFigFont{10}{14.4}{rm}$
	\lambda$}}}}} 
\put(1737,39){\makebox(0,0)[lb]{\smash{{{\SetFigFont{10}{14.4}{rm}$
	\sqrt{\lambda}$}}}}}
\end{picture}
}
$$
\caption{Case 1: $\lambda\neq 0, 1, \infty$.}
\end{figure}

{\sc Case 2:} $\lambda_s=0, j_s \neq 0, \infty$. 

Consider $\overline{C}$. It
has 4 sections, only 2 of which reduce to the origin and 2 reduce to another
point. The associated stable 4-pointed curve of genus 1, denoted  $C_s$, has
rational components $C_0^1$ and $C_0^2$ attached at these two points, through
which the sections 
pass. These rational components are branced over the rational components of
$D_s$ at the 2 sections. Moreover, the two attaching points  of  $\overline{C}$
are fixed points of the involution, therefore they must  be branch points of
$f_s$ on $C_0^1$ and $C_0^2$. This means that $C_0^i$ are ramified at 3 points,
which implies that they map purely inseparably onto the components of $D$.  The
component $C_1$ maps to the node of $D$.
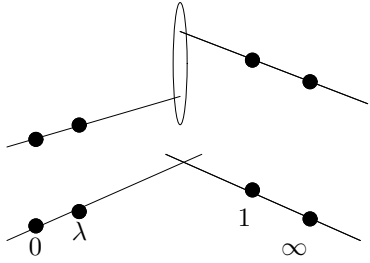
\begin{figure}[tbh]
$$
\setlength{\unitlength}{0.0005in}
\begingroup\makeatletter\ifx\SetFigFont\undefined
% extract first six characters in \fmtname
\def\x#1#2#3#4#5#6#7\relax{\def\x{#1#2#3#4#5#6}}%
\expandafter\x\fmtname xxxxxx\relax \def\y{splain}%
\ifx\x\y   % LaTeX or SliTeX?
\gdef\SetFigFont#1#2#3{%
  \ifnum #1<17\tiny\else \ifnum #1<20\small\else
  \ifnum #1<24\normalsize\else \ifnum #1<29\large\else
  \ifnum #1<34\Large\else \ifnum #1<41\LARGE\else
     \huge\fi\fi\fi\fi\fi\fi
  \csname #3\endcsname}%
\else
\gdef\SetFigFont#1#2#3{\begingroup
  \count@#1\relax \ifnum 25<\count@\count@25\fi
  \def\x{\endgroup\@setsize\SetFigFont{#2pt}}%
  \expandafter\x
    \csname \romannumeral\the\count@ pt\expandafter\endcsname
    \csname @\romannumeral\the\count@ pt\endcsname
  \csname #3\endcsname}%
\fi
\fi\endgroup
{\renewcommand{\dashlinestretch}{30}
\begin{picture}(3774,2687)(0,-10)
\put(1812,2027){\ellipse{150}{1274}}
\put(312,339){\blacken\ellipse{150}{150}}
\put(312,339){\ellipse{150}{150}}
\put(762,489){\blacken\ellipse{150}{150}}
\put(762,489){\ellipse{150}{150}}
\put(312,1239){\blacken\ellipse{150}{150}}
\put(312,1239){\ellipse{150}{150}}
\put(762,1389){\blacken\ellipse{150}{150}}
\put(762,1389){\ellipse{150}{150}}
\put(2562,2064){\blacken\ellipse{150}{150}}
\put(2562,2064){\ellipse{150}{150}}
\put(3162,1839){\blacken\ellipse{150}{150}}
\put(3162,1839){\ellipse{150}{150}}
\put(2562,714){\blacken\ellipse{150}{150}}
\put(2562,714){\ellipse{150}{150}}
\put(3162,414){\blacken\ellipse{150}{150}}
\put(3162,414){\ellipse{150}{150}}
\path(12,1164)(1812,1689)
\path(12,189)(2037,1089)
\path(1812,2364)(3762,1614)
\path(1812,2364)(3762,1614)
\path(1662,1089)(3687,189)
\path(1662,1089)(3687,189)
\put(237,39){\makebox(0,0)[lb]{\smash{{{\SetFigFont{10}{14.4}{rm}$0$}}}}}
\put(687,189){\makebox(0,0)[lb]{\smash{{{\SetFigFont{10}{14.4}{rm}$
	\lambda$}}}}}
\put(2412,339){\makebox(0,0)[lb]{\smash{{{\SetFigFont{10}{14.4}{rm}$1$}}}}}
\put(2862,39){\makebox(0,0)[lb]{\smash{{{\SetFigFont{10}{14.4}{rm}$
	\infty$}}}}}
\end{picture}
}
$$
\caption{Case 2: $\lambda_s=0, j_s \neq 0, \infty$.}
\end{figure}

{\sc Case 3:} $\lambda_s=0, j_s = \infty$. Considering   $\overline{C}$, we can
take 2 of the sections, one passing through the origin and one through the
node. The stabilization as a 2-pointed curve is a ``banana curve'' of two
rational components 
meeting in two points, with one marking on  each component. Adding the two
other sections and arguing as in case 2 we arrive at the following picture: 

The curve $C_s$ has 4 components, denoted $C_1^1$, $C_1^2$, $C_0^1$ and
$C_0^2$. The components $C_1^i$ are attached to each other at 2 points, and
$f_s$ maps them to the node of $D_s$. The component $C_1^i$ is also atatched
at one point to $C_0^i$. The morphism $f_s$ maps $C_0^i$ purely inseparably to
the components of $D$. 
\begin{figure}[tbh]
$$
\setlength{\unitlength}{0.0005in}
\begingroup\makeatletter\ifx\SetFigFont\undefined
% extract first six characters in \fmtname
\def\x#1#2#3#4#5#6#7\relax{\def\x{#1#2#3#4#5#6}}%
\expandafter\x\fmtname xxxxxx\relax \def\y{splain}%
\ifx\x\y   % LaTeX or SliTeX?
\gdef\SetFigFont#1#2#3{%
  \ifnum #1<17\tiny\else \ifnum #1<20\small\else
  \ifnum #1<24\normalsize\else \ifnum #1<29\large\else
  \ifnum #1<34\Large\else \ifnum #1<41\LARGE\else
     \huge\fi\fi\fi\fi\fi\fi
  \csname #3\endcsname}%
\else
\gdef\SetFigFont#1#2#3{\begingroup
  \count@#1\relax \ifnum 25<\count@\count@25\fi
  \def\x{\endgroup\@setsize\SetFigFont{#2pt}}%
  \expandafter\x
    \csname \romannumeral\the\count@ pt\expandafter\endcsname
    \csname @\romannumeral\the\count@ pt\endcsname
  \csname #3\endcsname}%
\fi
\fi\endgroup
{\renewcommand{\dashlinestretch}{30}
\begin{picture}(3774,2912)(0,-10)
\put(2393.250,2045.250){\arc{1762.899}{1.9869}{4.2012}}
\put(1168.731,2043.283){\arc{1886.993}{5.1717}{7.3038}}
\put(312,339){\blacken\ellipse{150}{150}}
\put(312,339){\ellipse{150}{150}}
\put(762,489){\blacken\ellipse{150}{150}}
\put(762,489){\ellipse{150}{150}}
\put(312,1239){\blacken\ellipse{150}{150}}
\put(312,1239){\ellipse{150}{150}}
\put(762,1389){\blacken\ellipse{150}{150}}
\put(762,1389){\ellipse{150}{150}}
\put(2562,2064){\blacken\ellipse{150}{150}}
\put(2562,2064){\ellipse{150}{150}}
\put(3162,1839){\blacken\ellipse{150}{150}}
\put(3162,1839){\ellipse{150}{150}}
\put(2562,714){\blacken\ellipse{150}{150}}
\put(2562,714){\ellipse{150}{150}}
\put(3162,414){\blacken\ellipse{150}{150}}
\put(3162,414){\ellipse{150}{150}}
\path(12,1164)(1887,1764)
\path(12,189)(2037,1089)
\path(1737,2439)(3762,1614)
\path(1737,2439)(3762,1614)
\path(1662,1089)(3687,189)
\path(1662,1089)(3687,189)
\put(237,39){\makebox(0,0)[lb]{\smash{{{\SetFigFont{10}{14.4}{rm}$0$}}}}}
\put(687,189){\makebox(0,0)[lb]{\smash{{{\SetFigFont{10}{14.4}{rm}$
	\lambda$}}}}}
\put(2412,339){\makebox(0,0)[lb]{\smash{{{\SetFigFont{10}{14.4}{rm}$1$}}}}}
\put(2862,39){\makebox(0,0)[lb]{\smash{{{\SetFigFont{10}{14.4}{rm}$
	\infty$}}}}}
\end{picture}
}
$$
\caption{Case 3: $\lambda_s=0, j_s = \infty$.}
\end{figure}
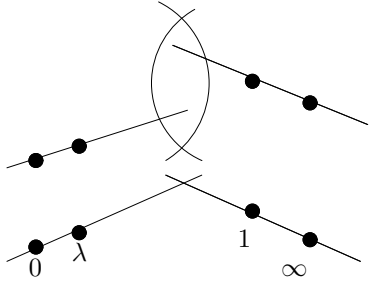

{\sc Case 4:}  $\lambda_s=0, j_s = 0$. The Legendre equation $y^2 =
x(x-1)(x-\lambda)$ is singular at $x=0$. Applying the transformation $x \mapsto
x/\lambda$ we get the Legendre equation $y^2 =
x(x-1)(x-\lambda^{-1})$ which  is singular at $x=\infty$. It follows that the
double cover of the stable 4-pointed curve $D_s$ is singular at the
node, and as in case 1 the singularity contributes genus 1. The stable limit
$C_s$ thus has to contain a component of genus 1 and $j$ invariant 0, attached
at one point to the rest of the curve. The only stable configuration possible
is the following: 

The curve $C_s$ has 4 components, $C_1$, $C_0'$, $C_0^1$ and 
$C_0^2$. The component $C_1$ has genus 1 and $j$-invariant 0. The other
components have genus 0. The component $C'$ is attached at one point to each of
$C_1$,  $C_0^1$ and 
$C_0^2$. Both  $C_1$ and $C_0'$ are mapped to the node of $D_s$, and $C_0^1$
and  
$C_0^2$ map purely inseparably to the two components of $D_s.$

\begin{figure}[tbh]
$$
\setlength{\unitlength}{0.0005in}
\begingroup\makeatletter\ifx\SetFigFont\undefined
% extract first six characters in \fmtname
\def\x#1#2#3#4#5#6#7\relax{\def\x{#1#2#3#4#5#6}}%
\expandafter\x\fmtname xxxxxx\relax \def\y{splain}%
\ifx\x\y   % LaTeX or SliTeX?
\gdef\SetFigFont#1#2#3{%
  \ifnum #1<17\tiny\else \ifnum #1<20\small\else
  \ifnum #1<24\normalsize\else \ifnum #1<29\large\else
  \ifnum #1<34\Large\else \ifnum #1<41\LARGE\else
     \huge\fi\fi\fi\fi\fi\fi
  \csname #3\endcsname}%
\else
\gdef\SetFigFont#1#2#3{\begingroup
  \count@#1\relax \ifnum 25<\count@\count@25\fi
  \def\x{\endgroup\@setsize\SetFigFont{#2pt}}%
  \expandafter\x
    \csname \romannumeral\the\count@ pt\expandafter\endcsname
    \csname @\romannumeral\the\count@ pt\endcsname
  \csname #3\endcsname}%
\fi
\fi\endgroup
{\renewcommand{\dashlinestretch}{30}
\begin{picture}(3774,3361)(0,-10)
\put(312,339){\blacken\ellipse{150}{150}}
\put(312,339){\ellipse{150}{150}}
\put(762,489){\blacken\ellipse{150}{150}}
\put(762,489){\ellipse{150}{150}}
\put(312,1239){\blacken\ellipse{150}{150}}
\put(312,1239){\ellipse{150}{150}}
\put(762,1389){\blacken\ellipse{150}{150}}
\put(762,1389){\ellipse{150}{150}}
\put(2562,2064){\blacken\ellipse{150}{150}}
\put(2562,2064){\ellipse{150}{150}}
\put(3162,1839){\blacken\ellipse{150}{150}}
\put(3162,1839){\ellipse{150}{150}}
\put(2562,714){\blacken\ellipse{150}{150}}
\put(2562,714){\ellipse{150}{150}}
\put(3162,414){\blacken\ellipse{150}{150}}
\put(3162,414){\ellipse{150}{150}}
\put(1812,3002){\ellipse{150}{676}}
\path(12,1164)(1887,1764)
\path(12,189)(2037,1089)
\path(1737,2439)(3762,1614)
\path(1737,2439)(3762,1614)
\path(1662,1089)(3687,189)
\path(1662,1089)(3687,189)
\path(1812,1539)(1812,2814)
\put(237,39){\makebox(0,0)[lb]{\smash{{{\SetFigFont{10}{14.4}{rm}$0$}}}}}
\put(687,189){\makebox(0,0)[lb]{\smash{{{\SetFigFont{10}{14.4}{rm}$
	\lambda$}}}}}
\put(2412,339){\makebox(0,0)[lb]{\smash{{{\SetFigFont{10}{14.4}{rm}$1$}}}}}
\put(2862,39){\makebox(0,0)[lb]{\smash{{{\SetFigFont{10}{14.4}{rm}$
	\infty$}}}}}
\end{picture}
}
$$
\caption{Case 4: $\lambda_s=0, j_s = 0$.}
\end{figure}
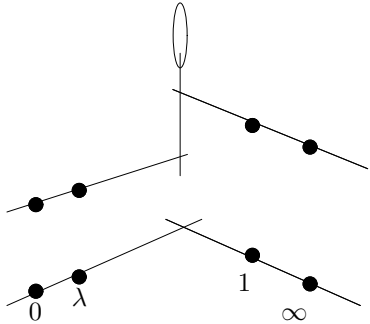

A situation identical to cases 1,3,4 holds for $\lambda_s=1,\infty$, by
permuting the marked points. 

\subsection{Speculations about intermediate quotient curves.}   
One would like to insert an intermediate curve $P$ between $C$ and $D$ so
that  $C \to P$ is finite. In the present case of $\ocH_{2,4}$ there is a
natural involution $\sigma$ on $C$, and in characteristic 0 we have $C /
\langle\sigma\rangle = D$. Thus in general, given a family $C \to D$ over a
scheme $T$ which is an element of $\ocH_{2,4}(T)$,  it is natural to consider
the quotient curve $P =C 
/ \langle\sigma\rangle$, and we have natural morphisms 
$C \to P \to D$. 

 Unfortunately {\em the formation of $P$ does not commute
 with  base change!} Indeed, if $T$ has generic point of characteristic 0 and
special point in characteristic 2 as before, then over the generic points of
$D_s$, the morphism $P \to D$ is an isomorphism. However if $T$ is replaced by
$s$, then $\sigma$ acts {\em trivially} on the components $C_0^i$, therefore $P
\to D$ inseparable of degree 2  over the generic points of
$D$! 

There is, however, a remedy. In analogy to the methods of \cite{av}, one could
take $P$ to be the {\em stack}  $C / \langle\sigma\rangle$. This stack
automatically commutes with base changes. The morphism $C \to P$ is a principal
$\ZZ/2\ZZ$ bundle, giving rise to a representable morphism to the classifying
stack $P \to B(\ZZ/2\ZZ)$. In a sense yet to be understood, this may be
considered a stable map.

In order to introduce intermediate curves for general complete Hurwitz stacks,
one might replace a Hurwitz cover $C \to D$ by its associated $S_d$-cover
$\tilde{C} \to D$,
where $S_d$ is the symmetric group on $d$ letters. This again is in analogy
with \cite{av}. One can complete the locus of these $S_d$-covers inside the
appropriate stack of $S_d$-equivariant stable maps, and obtain a new stack
which one might denote $\ocH_{S_d, h, n}$. An object of this stack over a
reduced base 
scheme $T$ is an
$S_n$-equivariant $T$-morphism $\tilde{C} \to D$, branched over $n$-marked
points 
in $D$, whose geometric fibers can be lifted to Horwitz covers in
characteristic 0. One can 
obtain a replacement 
for $C$ by taking the quotient  $\tilde{C}/S_{d-1}$. This quotient should be
taken as a Deligne-Mumford stack. The quotient
$\tilde{C}/S_{d}$ is a candidate for $P$. 

%------------------------------------------------------------

\end{document}